\documentclass[1p]{elsarticle}
\usepackage{amssymb}
\usepackage{amsfonts}

\newtheorem{theorem}{Theorem}
\newtheorem{conjecture}[theorem]{Conjecture}
\newtheorem{corollary}[theorem]{Corollary}

\newtheorem{observation}[theorem]{Observation}
\newtheorem{claim}
{Claim}

\newproof{pf}{Proof}

\begin{document}
\title{Distant set distinguishing edge colourings of graphs}

\author{Jakub Przyby{\l}o\fnref{fn1,fn2}}
\ead{jakubprz@agh.edu.pl, phone: 048-12-617-46-38,  fax: 048-12-617-31-65}

\fntext[fn1]{Financed within the program of the Polish Minister of Science and Higher Education
named ``Iuventus Plus'' in years 2015-2017, project no. IP2014 038873.}
\fntext[fn2]{Partly supported by the Polish Ministry of Science and Higher Education.}

\address{AGH University of Science and Technology, al. A. Mickiewicza 30, 30-059 Krakow, Poland}

\begin{abstract}
We consider the following extension of the concept of adjacent strong edge colourings
of graphs without isolated edges.
Two distinct vertices which are at distant at most $r$ in a graph are called \emph{$r$-adjacent}.
The least number of colours in a proper edge colouring
of a graph $G$ such that the sets of colours met by any $r$-adjacent vertices in $G$ are distinct
is called the \emph{$r$-adjacent strong chromatic index} of $G$ and denoted by $\chi'_{a,r}(G)$.
It has been conjectured that $\chi'_{a,1}(G)\leq\Delta+2$
if $G$ is connected of maximum degree $\Delta$ and non-isomorphic to $C_5$,
while Hatami proved that there is a constant $C$, $C\leq 300$, such that $\chi'_{a,1}(G)\leq\Delta+C$
if $\Delta>10^{20}$ [J. Combin. Theory Ser. B 95 (2005) 246--256].
We conjecture that a similar statement should hold for any $r$, i.e., that
for each positive integer $r$ there exist constants $\delta_0$ and $C$ such that
$\chi'_{a,r}(G) \leq \Delta+C$ for every graph without an isolated edge and with minimum degree $\delta \geq \delta_0$,
and argue that a lower bound on $\delta$ is unavoidable in such a case (for $r>2$).
Using the probabilistic method we prove such upper bound to hold for graphs with $\delta\geq \epsilon\Delta$,
for every $r$ and any fixed $\varepsilon\in(0,1]$, i.e., in particular for regular graphs.
We also support the conjecture by proving an upper bound $\chi'_{a,r}(G) \leq (1+o(1))\Delta$
for graphs with $\delta\geq r+2$. 
\end{abstract}


\maketitle

\section{Introduction}
In~\cite{Zhang} Zhang et al. posed the following problem.
Consider a graph $G=(V,E)$ containing no isolated edges and its proper edge colouring $c:E\to\{1,2,\ldots,k\}$.
For any $v\in V$ denote by $S_c(v)$ the set of colours incident with $v$,
i.e.,
$$S_c(v):=\{c(uv):u\in N(v)\},$$ where $N(v)$ is the set of neighbours of $v$.
We shall also denote $S_c(v)$ simply by $S(v)$
when this causes no ambiguities, and
refer to it as the \emph{colour pallet} of $v$.
We call vertices $u$ and $v$ \emph{distinguished} if $S(u)\neq S(v)$.
The least number of colours in a proper edge colouring $c$
which distinguishes the ends of all edges of $G$,
i.e. such that $S_c(u)\neq S_c(v)$ for every $uv\in E$,
is called the \emph{neighbour set distinguishing index} or the \emph{adjacent strong chromatic index} and denoted by $\chi'_a(G)$,
see~\cite{Akbari,BalGLS,FlandrinMPSW,Hatami,Zhang}
(also for other notations used).
\begin{conjecture}[\cite{Zhang}]\label{ZhangsConjecture}
For every connected graph $G$, $\chi'_{a}(G)\leq \Delta(G)+2$, unless $G$ is isomorphic to $K_2$ or $C_5$.
\end{conjecture}
This problem remains open despite many articles studying it.
In particular, it has been showed that $\chi'_a(G)\leq 3\Delta(G)$, \cite{Akbari}, and $\chi'_a(G)\leq \Delta(G)+O(\log\chi(G))$, \cite{BalGLS}.
Moreover, the conjecture was verified to hold for special families of graphs, like e.g.
 bipartite graphs or graphs of maximum degree 3, see~\cite{BalGLS}.
The following thus far best general upper bound is due to Hatami.
\begin{theorem}[\cite{Hatami}]
If $G$ is a graph with no isolated edges and maximum degree $\Delta>10^{20}$, then $\chi'_{a}(G)\leq \Delta+300$.
\end{theorem}
Note that this implies that $\chi'_{a}(G)\leq \Delta(G)+C$ for every graph $G$ containing no isolated edges, where $C$ is some constant.

Let $r$ be any positive integer. Vertices $u,v$ of $G$ shall be called \emph{$r$-neighbours} (or \emph{$r$-adjacent})
if $1\leq d(u,v)\leq r$, where $d(u,v)$ denotes the distance of $u$ and $v$ in $G$.
In this paper, similarly as e.g. within the concept of distant chromatic numbers (see~\cite{DistChrSurvey} for a survey of this topic),
we propose an extension of the study above towards distinguishing
not only neighbours, but also vertices at some limited distance (from each other).
The least number of colours in a proper edge colouring $c$ of $G$ such that $S_c(u)\neq S_c(v)$
for every pair of vertices $u,v\in V$ with $1\leq d(u,v)\leq r$,
so-called \emph{$r$-distant set distinguishing colouring}
(or \emph{$r$-adjacent strong edge colouring}),
shall be called the \emph{$r$-distant set distinguishing index} or
\emph{$r$-adjacent strong chromatic index}, and denoted by $\chi'_{a,r}(G)$.
This graph invariant has already been considered in \cite{KemnitzMarangio} and~\cite{MocSotak} under the name of \emph{$d$-strong chromatic index}
(see \cite{Akbari,Tian,ZhangDistant} for other notations),
mainly 
for paths, cycles and circulant graphs,
aiming towards providing
a series of counterexamples to a conjecture from~\cite{ZhangDistant}.
%
Some aspects of this concept were also 
investigated in~\cite{Akbari}
with respect to trees and small values of $r$.
As for general upper bounds, in~\cite{Tian} it was proved that
$\chi'_{a,2}(G)\leq 32(\Delta(G))^2$ if $\Delta(G)\geq 4$,
$\chi'_{a,3}(G)\leq 8(\Delta(G))^\frac{5}{2}$ if $\Delta(G)\geq 6$
and $\chi'_{a,r}(G)\leq 2\sqrt{2(r-1)}(\Delta(G))^\frac{r+2}{2}$
for $r\geq 4$ if $\Delta(G)\geq 4$.
Within this paper, among others, we intend to improve these bounds significantly under some (unavoidable) degree conditions.

The cornerstone of the general field of vertex distinguishing colourings,
which is rich in many interesting open problems and conjectures, is the graph invariant called the irregularity strength.
Consider a (not necessarily proper) edge colouring $c:E\to\{1,2,\ldots,k\}$
of a graph $G=(V,E)$ containing no isolated edges. Denote by
$$s_c(v)=\sum_{u\in N(v)}c(uv)$$
the sum of colours incident with any vertex $v\in V$.
The \emph{irregularity strength} of $G$, $s(G)$, is then the least integer $k$ admitting
such $c$ with $s_c(u)\neq s_c(v)$ for all $u,v\in V$, $u\neq v$, see~\cite{Chartrand}.
Note that equivalently, it is equal to the minimal integer $k$ so that we are able
to construct an irregular multigraph (a multigraph with pairwise distinct degrees of all the vertices)
of $G$ by multiplying some of its edges, each at most $k$ times.
Out of the extensive bibliography devoted to this parameter,
it is in particular worth mentioning~\cite{Aigner} and~\cite{Nierhoff}, where the sharp upper bound $s(G)\leq n-1$ (with $n=|V|$) was settled,
and~\cite{KalKarPf}, where a better upper bound $s(G)\leq 6\lceil n/\delta\rceil$ for graphs with sufficiently large minimum degree $\delta$ is proved. See~\cite{Lehel} for an interesting survey and open problems in
this topic as well.
This problem
is also related with the study of irregular graphs by
Chartrand, Erd\H{o}s and Oellermann, see~\cite{ChartrandErdosOellermann},
and gave rise to many other intriguing graph invariants.
In particular the following concept
was closely related to the later local version of irregularity strength
--
the well known problem commonly referred to as 1--2--3--Conjecture, see~\cite{123KLT}
by Karo\'nski, \L uczak and Thomason (and~\cite{KalKarPf_123} for the best result concerning this).
Let $K$ be the least integer $K$ so that for every graph $G=(V,E)$ without isolated edges there
exists a (not necessarily proper) edge colouring
$c:E\to\{1,2,\ldots,K\}$
such that for each edge $uv\in E$, the multisets of colours incident with $u$ and $v$ are distinct.
Note that this problem is a natural correspondent of the concept of $\chi'_a(G)$,
as the only difference is the requirement concerning the properness of the colourings investigated.
However, by~\cite{Aldred} it is known that $K\leq 4$, or even $K\leq 3$ suffices for graphs with minimum degree $\delta\geq 10^3$,
while $\chi'_{a}$ is bounded from below by $\Delta$.
In fact one of our main motivations for
studying the parameters $\chi'_{a,r}$ in this paper
is a desire to expose the leading impact
of the required properness of colourings on
the number of colours needed to distinguish
vertices by their incident (multi)sets, and that
usually not many more colours are needed
if we wish to distinguish not only neighbours,
but also vertices at distance $2$, $3$,..., $r<\infty$.
We believe, and
prove in many cases that for every fixed $r$, if only $\delta(G)$ is `not very small', then
$\chi'_{a,r}(G)\leq \Delta+C$ for each graph $G$ without an isolated edge,
where $C$ is some constant dependent on $r$, cf. Conjecture~\ref{Przybylo_main_conjecture} and Theorem~\ref{przybylo_main_bound} below.
We also confirm our conjecture asymptotically by proving an upper bound $\chi'_{a,r}(G)\leq (1+o(1))\Delta$, cf. Theorem~\ref{first_attempt_theorem}.

Moreover, or maybe even more importantly, by our research, we wish to reveal a
difference between the two main concepts of distinguishing
vertices (at a bounded distance) of the field, i.e.,
this with respect to sets and this based on sums.
Intriguingly this expected difference was elusive while
distinguishing only neighbours.
The least integer $k$ so that a proper edge colouring $c:E\to\{1,2,\ldots,k\}$
exists with $s_c(u)\neq s_c(v)$ for every edge $uv\in E$ is denoted by $\chi'_{\sum}(G)$.
Note that $\chi'_a(G)\leq\chi'_{\sum}(G)$ for every graph $G$ without isolated edges.
Though the requirement $s_c(u)\neq s_c(v)$ is much stronger than $S_c(u)\neq S_c(v)$,
it was conjectured in~\cite{FlandrinMPSW} that $\chi'_{\sum}(G)\leq \Delta(G)+2$
for every connected graph non-isomorphic to $K_2$ nor $C_5$ (similarly as for $\chi'_a$ in~\cite{Zhang}, cf. Conjecture~\ref{ZhangsConjecture}),
what was also asymptotically confirmed in~\cite{Przybylo_asym_optim}, where it was proved that $\chi'_{\sum}(G)\leq (1+o(1))\Delta$.
In fact (almost) all exact values of the both parameters, settled for some special families of graphs coincide,
see~\cite{FlandrinMPSW} for further comments.
See also~\cite{BonamyPrzybylo,FlandrinMPSW,Przybylo_CN_1,Przybylo_CN_2} for other results concerning $\chi'_{\sum}$.
On the other hand, if we consider a distant version of the problem in sum environment,
even setting aside the properness of colourings,
then the number of colours required grows rapidly with $r$.
For any positive integer $r$, let $s_r(G)$ be the least integer $k$ so that an edge colouring $c:E\to\{1,2,\ldots,k\}$
exists with $s_c(u)\neq s_c(v)$ for every pair of $r$-neighbours $u,v$.
Then it is known that $s_r(G)\leq 6\Delta^{r-1}$ for every graph $G$ without isolated edges, see~\cite{Przybylo_distant}.
On the other hand, it can be proved that there are graphs for which
the parameter $s_r(G)$ cannot be much smaller than $\Delta^{r-1}$
for arbitrarily large $\Delta$, also in the family of regular graphs,
contrary to the parameter $\chi'_{a,r}$ --
cf. Conjecture~\ref{Przybylo_main_conjecture}, Theorems~\ref{first_attempt_theorem}, \ref{przybylo_main_bound}, and
especially Corollary~\ref{regular_graphs_corollary} below.

\section{Lower Bound and Main Conjecture}

We shall use the following well known inequalities:
\begin{equation}\label{Choose_inequalities}
\left(\frac{a}{b}\right)^b\leq{a\choose b}
\leq\frac{a^b}{b!}\leq \left(\frac{ea}{b}\right)^b.
\end{equation}

Suppose we wish to prove that $\chi'_{a,r}(G)\leq \Delta+C$ (or at least $\chi'_{a,r}(G)\leq (1+o(1))\Delta$),
then some (minor) assumptions concerning the structure of $G$, e.g. the minimum degree of $G$, are unavoidable.
\begin{observation}
Let $r>2$ be an integer. Suppose $\chi'_{a,r}(G)\leq \Delta(G)+C$ (or at least $\chi'_{a,r}(G)\leq (1+o(1))\Delta$)
for every graph $G$ with $\delta(G)\geq \delta_0$ and without isolated edges, where $C$ is some constant.
Then we must have $\delta_0\geq r$ if $r\geq 7$, or at least $\delta_0\geq r-1$ otherwise.
\end{observation}
\begin{pf}
Assume first that $r\geq 9$.
In order to prove the thesis we construct an infinite family of graphs with minimum degree $\delta = r-1$
for which there exist no $r$-distant set distinguishing colourings with just $\Delta+C$ (or $(1+o(1))\Delta$) colours.
We shall use so called \emph{undirected de Bruijn graph} of type $(t,k)$, whose vertex set is formed by all sequences of length $k$
the entries of which are taken from a fixed alphabet consisting of $t$ distinct letters, and in which two distinct vertices $(a_1,\ldots,a_k)$ and $(b_1,\ldots,b_k)$ are joined by an edge if either $a_i = b_{i+1}$ for $1\leq i \leq k - 1$, or if $a_{i+1} = b_i$ for $1 \leq i \leq k - 1$ (or $k=1$).
Such a graph, which we shall denote by $D_{t,k}$,
has maximum degree $\Delta(D_{t,k})\leq 2t$, 
order $t^k$ and diameter $k$, 
and provides a nontrivial lower bound in the study of so called \emph{Moore bound}, concerning the largest order of a graph with given maximum degree and diameter, see e.g. a survey by Miller and {\v S}ir{\'a}{\v n}~\cite{Mirka}.
For any positive integer $N$ (which shall be required to be large enough
later on), consider the graph $G'$ obtained
by taking $2N[N(r-2)(r-1)]^{r-2}$ disjoint copies of $K_{r}$ and
identifying exactly one vertex from each of these with some vertex of $D_{N(r-2)(r-1),r-2}$ (of order $[N(r-2)(r-1)]^{r-2}$) so that each vertex from $D_{N(r-2)(r-1),r-2}$ is incident with exactly $2N$ such complete graphs.
Note that $\Delta(G')\leq 2N(r-2)(r-1)+2N(r-1)=2N(r-1)^2$ and $G'$ contains $[N(r-2)(r-1)]^{r-2}2N(r-1)=2[N(r-1)]^{r-1}(r-2)^{r-2}$
vertices of (minimum) degree $r-1$, every two of which are at distance at most $r$.
On the other hand, with only $\Delta(G')+C$ colours admitted,
any given vertex of degree $r-1$ may be assigned one of at most ${\Delta(G')+C \choose r-1} \leq {2N(r-1)^2+C \choose r-1}$
potential colour pallets (where $C$ is some fixed constant).
To prove that there exists no $r$-distant set distinguishing colouring of this graph with $\Delta(G')+C$ colours
it is thus sufficient to show that
\begin{equation}\label{delta_lower_bound_ineq}
2[N(r-1)]^{r-1}(r-2)^{r-2} > {2N(r-1)^2+C \choose r-1}.
\end{equation}
Since by~(\ref{Choose_inequalities}),
\begin{equation}\label{delta_lower_bound_ineq_2}
{2N(r-1)^2+C \choose r-1}\leq\left(\frac{e[2N(r-1)^2+C]}{r-1}\right)^{r-1}\leq [5.5N(r-1)]^{r-1}
\end{equation}
(for $N$ sufficiently large),
by (\ref{delta_lower_bound_ineq}) and (\ref{delta_lower_bound_ineq_2}) it is then sufficient to prove the inequality
$$2[N(r-1)]^{r-1}(r-2)^{r-2} > [5.5N(r-1)]^{r-1},$$
or equivalently
$$\left(\frac{r-2}{5.5}\right)^{r-2} > 2.75,$$
which holds for $r\geq 9$ (the left hand side above is an increasing function of $r$ for $r\geq 9$).

Note that the same argument holds even if we admit $\Delta(G')(1+o(1))$ instead of $\Delta(G')+C$ colours (as inequality (\ref{delta_lower_bound_ineq_2}) holds also after such substitution).
Moreover, by more careful estimations, one can easily show that in both cases the same is also true already for $r\geq 7$ (by proving directly that~(\ref{delta_lower_bound_ineq}) holds for $r=7,8$ and $N$ sufficiently large).

Finally, for $r\leq 6$, by a similar approach, one can show that we need at least the assumption $\delta\geq r-1$.
It is sufficient to use complete graphs $K_{r-1}$, instead of $K_r$ in the construction above.
Then for any fixed $r\in\{3,4,5,6\}$, $G'$
would contain $\Omega(\Delta^{r-1})$
vertices of degree $r-2$ every pair of which would be at distance at most $r$,
while there would be $O(\Delta^{r-2})$ available colour pallets for these vertices (for $N$ tending to infinity).
$\blacksquare$
\end{pf}

We conclude this section by posing the following general conjecture, where we believe that $\delta_0$
should be roughly equal to $r$ (up to some small additive constant).
\begin{conjecture}\label{Przybylo_main_conjecture}
For each positive integer $r$ there exist constants $\delta_0$ and $C$ such that
$$\chi'_{a,r}(G) \leq \Delta(G)+C$$
for every graph without an isolated edge and with $\delta(G) \geq \delta_0$.
\end{conjecture}

\section{Asymptotic Confirmation}


First we shall prove that $(1+o(1))\Delta$ colours are sufficient in case of
graphs with minimum degree larger than $r+1$. 
The proof is based on
the Lov\'asz Local Lemma, see e.g.~\cite{AlonSpencer}, combined with the Chernoff Bound, see e.g.~\cite{MolloyReed}.
\begin{theorem}[\textbf{The Local Lemma}]
\label{LLL-symmetric}
Let $A_1,A_2,\ldots,A_n$ be events in an arbitrary pro\-ba\-bi\-li\-ty space.
Suppose that each event $A_i$ is mutually independent of a set of all the other
events $A_j$ but at most $D$, and that ${\rm \emph{\textbf{Pr}}}(A_i)\leq p$ for all $1\leq i \leq n$. If
$$ ep(D+1) \leq 1,$$
then $ {\rm \emph{\textbf{Pr}}}\left(\bigcap_{i=1}^n\overline{A_i}\right)>0$.
\end{theorem}
\begin{theorem}[\textbf{Chernoff Bound}]\label{ChernofBoundTh}
For any $0\leq t\leq np$:
$${\rm\emph{\textbf{Pr}}}(|{\rm BIN}(n,p)-np|>t)<2e^{-\frac{t^2}{3np}},$$
where ${\rm BIN}(n,p)$ is the sum of $n$ independent variables, each equal to $1$ with probability $p$ and $0$ otherwise.
\end{theorem}

A colouring $c'$ which assigns colours to some part of the edges of a graph $G$ shall be called a \emph{partial colouring}.
Given such partial colouring $c'$ and a vertex $v$ of $G$, the set
$S_{c'}(v)$ is defined the same as for a complete edge colouring,
but only coloured edges are taken into account in this set.
By $d_{c'}(v)$ we shall also denote the number of edges incident with $v$ which are coloured under $c'$.

\begin{theorem}\label{first_attempt_theorem}
For every $r\geq 2$ and every graph $G$ of maximum degree $\Delta$ with 
$\delta(G)\geq r+2$ and without isolated edges,
$$\chi'_{a,r}(G)\leq (1+o(1))\Delta.$$
\end{theorem}

\begin{pf}
Let us fix a positive integer $r$ and suppose $G=(V,E)$ is a graph without isolated edges,
with maximum and minimum degrees $\Delta$ and $\delta$, resp., such that 
$\delta\geq r+2$.
We shall in fact prove that if $\Delta$ is sufficiently large, then
$$\chi'_{a,r}(G)\leq \left\lceil\frac{\Delta}{\ln^4\Delta}\right\rceil \left(\lceil\ln^4\Delta\rceil+5\lceil\ln^3\Delta\rceil+2\right).$$
Whenever needed we shall thus assume that $\Delta$ is large enough.

We shall use colours $1,2,\ldots,\lceil\frac{\Delta}{\ln^4\Delta}\rceil (\lceil\ln^4\Delta\rceil+5\lceil\ln^3\Delta\rceil+2)$,
which we arbitrarily partition into $t=\lceil\frac{\Delta}{\ln^4\Delta}\rceil$ subsets $C_1,C_2,\ldots,C_t$ of equal cardinalities.
Each $C_i$ is furthermore partitioned into two subsets $C'_i$ and $C''_i$, where $|C'_i|=\lceil\ln^4\Delta\rceil+\lceil\ln^3\Delta\rceil+1$
and $|C''_i|=4\lceil\ln^3\Delta\rceil+1$ for $i=1,\ldots,t$.
First we use probabilistic approach to construct a not necessarily proper edge colouring
$q:E\to\{1,2,\ldots,t\}$ attributing every edge an index of a set from the partition above.
For this aim we randomly and equi-probably choose the value $q(e)\in\{1,2,\ldots,t\}$ independently for every edge $e\in E$.
This part of the construction is designed to assure distinction for vertices of small degrees.
Thus for every pair of $r$-neighbours $u,v$ with $d(u)=d(v)\leq \ln^3\Delta$ in $G$, let
$D_{\{u,v\}}$ denote the event that $S_q(u)=S_q(v)$ (where $S_q(u),S_q(v)$ are the sets defined in the same manner as for proper edge colourings above). Since $|S_q(u)|\leq \ln^3\Delta$ then, and $u$ and $v$ might have at most one common edge (and $d(v)\geq r+2$),
\begin{equation}\label{asymptotois_probabil_1}
\mathbf{Pr}(D_{\{u,v\}}) \leq \left(\frac{\ln^3\Delta}{\lceil\Delta\ln^{-4}\Delta\rceil}\right)^{r+1}
\leq \left(\frac{\ln^7\Delta}{\Delta}\right)^{r+1}.
\end{equation}

Set $Q_i:=\{e\in E: q(e)=i\}$ for $i=1,\ldots,t$.
For the sake of our construction, none of the index colours may appear too many times around any vertex.
Hence, for every vertex $v\in V$ and a colour $i\in\{1,\ldots,t\}$, denote by $D_{v,i}$
the event that 
more than $\ln^4\Delta+\ln^3\Delta$ edges 
incident with $v$ belong to $Q_i$ (i.e., these edges are coloured with $i$).
Note that then, by the Chernoff Bound:
\begin{eqnarray}
\mathbf{Pr}(D_{v,i}) & = &
{\rm\textbf{Pr}}\left({\rm BIN}\left(d(v),\frac{1}{\lceil\Delta\ln^{-4}\Delta\rceil}\right)>\ln^4\Delta+\ln^3\Delta\right)\nonumber\\
& \leq & {\rm\textbf{Pr}}\left({\rm BIN}\left(\Delta,\frac{\ln^4\Delta}{\Delta}\right)>\ln^4\Delta+\ln^3\Delta\right)\nonumber\\
& \leq & {\rm\textbf{Pr}}\left(\left|{\rm BIN}\left(\Delta,\frac{\ln^4\Delta}{\Delta}\right)-\ln^4\Delta\right|>\ln^3\Delta\right)\nonumber\\
& \leq & 2e^{-\frac{\ln^6\Delta}{3\ln^4\Delta}} \leq 
\left(\frac{\ln^7\Delta}{\Delta}\right)^{r+1}.\label{asymptotois_probabil_2}
\end{eqnarray}

On the other hand, every event $D_{v,i}$ is mutually independent of all but at most
$(\Delta+1)t+\Delta \ln^3\Delta\cdot\Delta^{r-1}\leq 2\Delta^r\ln^3\Delta$
other events of both types considered above,
while every event $D_{\{u,v\}}$ is mutually independent of all but at most
$2(\ln^3\Delta+1)t+2(\ln^3\Delta+1)\ln^3\Delta\cdot\Delta^{r-1}\leq 2\Delta^r\ln^3\Delta$
other events of both types. Hence, for $\Delta$ sufficiently large,
by (\ref{asymptotois_probabil_1}), (\ref{asymptotois_probabil_2}) and the Local Lemma (with positive probability),
we may choose the colouring $q$ so that:
\begin{itemize}
\item[(i)]
$S_q(u)\neq S_q(v)$ for every pair of $r$-neighbours $u,v$ with $d(u)=d(v)\leq \ln^3\Delta$ in $G$;
\item[(ii)] $|\{u\in N(v):q(uv)=i\}| \leq \ln^4\Delta+\ln^3\Delta$ for every $v\in V$ and $i\in\{1,\ldots,t\}$.
\end{itemize}

Since by (ii) above, for every $i$,
the subgraph induced in $G$ by the edges coloured with $i$ has maximum degree at most $\ln^4\Delta+\ln^3\Delta$,
by Vizing's Theorem, we may properly recolour the edges of this subgraph using colours from $C'_i$.
We denote by $c$ the proper edge colouring of the entire graph $G$ obtained after performing the above for every $i=1,\ldots,t$.
%
Now randomly and independently we uncolour every edge of $G$,
each with probability $\frac{2}{\ln\Delta}$. Denote by $c'$ the obtained partial colouring of $G$.
For every vertex $v$, we denote the set of uncoloured edges incident with $v$ by $U_{c'}(v)$.
Note that if $d(v)=d\geq \ln^3\Delta$, then
$\mathbf{E}(|U_{c'}(v)|)=\frac{2d}{\ln\Delta}$
and thus by the Chernoff Bound,
\begin{equation}\label{A_v_estimation1}
\mathbf{Pr}\left(\left||U_{c'}(v)|-\frac{2d}{\ln\Delta}\right|>\frac{d}{\ln\Delta}\right)<2e^{-\frac{d}{6\ln\Delta}}\leq 2e^{-\frac{\ln^2\Delta}{6}}\leq \frac{1}{\Delta^{r+3}}.
\end{equation}
Analogously, since by (ii) above, $\mathbf{E}(|U_{c'}(v)\cap Q_i|)\leq\frac{2(\ln^4\Delta+\ln^3\Delta)}{\ln\Delta} \leq 3\ln^3\Delta$,
by the Chernoff Bound,
\begin{equation}\label{A_v_estimation2}
\mathbf{Pr}\left(|U_{c'}(v)\cap Q_i|>4\ln^3\Delta\right)<2e^{-\frac{(\ln^3\Delta)^2}{9\ln^3\Delta}}\leq \frac{1}{\Delta^{r+3}}.
\end{equation}
For every vertex $v$ with $d(v)=d\geq \ln^3\Delta$,
denote by $A_{v,0}$ the event that $||U_{c'}(v)|-\frac{2d}{\ln\Delta}|>\frac{d}{\ln\Delta}$,
and for $i=1,\ldots,t$, denote by $A_{v,i}$ the event that $|U_{c'}(v)\cap Q_i|>4\ln^3\Delta$.
For any two distinct vertices $u,v$ at distance at most $r$ in $G$ which are of the same degree $d$ with $\ln^3\Delta\leq d\leq\Delta$,
denote by $B_{\{u,v\}}$ the event that $|U_{c'}(u)|,|U_{c'}(v)|\in \left[\frac{d}{\ln\Delta},\frac{3d}{\ln\Delta}\right]$
and $S_{c'}(u)=S_{c'}(v)$.
Then (as $u$ and $v$ have at most one common incident edge):
\begin{eqnarray}
\mathbf{Pr}(B_{\{u,v\}}) &\leq& \mathbf{Pr}\left(S_{c'}(u)=S_{c'}(v)\wedge |U_{c'}(v)|\in \left[\frac{d}{\ln\Delta},\frac{3d}{\ln\Delta}\right]\right)\nonumber\\
&\leq& \mathbf{Pr}\left(S_{c'}(u)=S_{c'}(v)\left||U_{c'}(v)|\in \left[\frac{d}{\ln\Delta},\frac{3d}{\ln\Delta}\right]\right.\right)\nonumber\\
&\leq&
\sum_{j =\lceil\frac{d}{\ln\Delta}\rceil}^{\lfloor\frac{3d}{\ln\Delta}\rfloor}
\mathbf{Pr}\left(S_{c'}(u)=S_{c'}(v)\left||U_{c'}(v)|=j\right.\right)
\mathbf{Pr}\left(|U_{c'}(v)|=j\left||U_{c'}(v)|\in \left[\frac{d}{\ln\Delta},\frac{3d}{\ln\Delta}\right]\right.\right)\nonumber\\
&\leq&
\sum_{j =\lceil\frac{d}{\ln\Delta}\rceil}^{\lfloor\frac{3d}{\ln\Delta}\rfloor}
\left(\frac{2}{\ln\Delta}\right)^{j-1}
\mathbf{Pr}\left(|U_{c'}(v)|=j\left||U_{c'}(v)|\in \left[\frac{d}{\ln\Delta},\frac{3d}{\ln\Delta}\right]\right.\right)\nonumber\\
&\leq&
\left(\frac{2}{\ln\Delta}\right)^{\frac{d}{\ln\Delta}-1}
\sum_{j =\lceil\frac{d}{\ln\Delta}\rceil}^{\lfloor\frac{3d}{\ln\Delta}\rfloor}
\mathbf{Pr}\left(|U_{c'}(v)|=j\left||U_{c'}(v)|\in \left[\frac{d}{\ln\Delta},\frac{3d}{\ln\Delta}\right]\right.\right)\nonumber\\
&\leq& \left(\frac{2}{\ln\Delta}\right)^{\ln^2\Delta-1}\cdot 1
\leq \left(\frac{1}{e}\right)^{(r+3)\ln\Delta}
=\frac{1}{\Delta^{r+3}}.\label{B_uv_estimation}
\end{eqnarray}
Note that every event
of the form $A_{v,i}$ ($i\in\{0,1,\ldots,t\}$) or $B_{\{u,v\}}$ is mutually independent of all other events
of these forms indexed by vertices each of which is at distance at least
$2$ from both $v$ and $u$ (in the case of $B_{\{u,v\}}$), i.e., of all but at most
$2\cdot(t+2)(\Delta+1) \Delta^r \leq \Delta^{r+2}$ other events.
Moreover, by~(\ref{A_v_estimation1}),~(\ref{A_v_estimation2}) and~(\ref{B_uv_estimation}) each of these events occurs with probability at most $\frac{1}{\Delta^{r+3}}$.
Therefore, by the Lov\'asz Local Lemma,
we may commit our uncolourings so that none of these events holds for the obtained $c'$.
Then, $S_{c'}(u)\neq S_{c'}(v)$ for every pair of
$r$-neighbours $u,v\in V$ with $d(u)=d(v)\geq \ln^3\Delta$
in $G$. Moreover, as $A_{v,i}$ does not hold for every $v\in V$ and $i\in\{1,\ldots,t\}$,
then the subgraph induced in $G$ by the uncoloured edges which belong to $Q_i$
has maximum degree at most $4\ln^3\Delta$, and thus can be coloured properly with
(yet unused) colours from $C''_i$.
After colouring each such subgraph for $i=1,\ldots,t$, we obtain our final proper edge colouring $c''$ of $G$.
%
Since we have used new colours, we still have $S_{c''}(u)\neq S_{c''}(v)$ whenever $u,v\in V$ are $r$-neighbours
one of which has degree at least $\ln^3\Delta$.
Otherwise, the same holds by the condition (i) above (as $c''(e)\in C_i$ for every edge $e$ with $q(e)=i$). $\blacksquare$
\end{pf}

\section{Almost Optimal Upper Bound}

In the following we shall prove an upper bound $\chi'_{a,r}(G)\leq\Delta+C$,
which is optimal up to the additive constant $C$,
for $\delta$ linear in $\Delta$, see Theorem~\ref{przybylo_main_bound} below.
For this aim we shall need one additional probabilistic tool, see e.g.~\cite{MolloyReed}.
\begin{theorem}[\textbf{Talagrand's Inequality}]
Let $X$ be a non-negative random variable, not identically $0$, which is determined by
$l$ independent trials $T_1,\ldots,T_l$, and satisfying the following for some $c,k>0$:
\begin{itemize}
\item[1.] changing the outcome of any one trial can affect $X$ by at most $c$, and
\item[2.] for any $s$, if $X\geq s$ then there is a set of at most $ks$ trials whose outcomes certify that $X\geq s$,
\end{itemize}
then for any $0\leq t\leq \mathbf{E}(X)$,
$${\rm\emph{\textbf{Pr}}}(|X-{\rm\emph{\textbf{E}}}(X)|>t+60c\sqrt{k{\rm\emph{\textbf{E}}}(X)})\leq 4e^{-\frac{t^2}{8c^2k\mathbf{E}(X)}}.$$
\end{theorem}

\begin{theorem}\label{przybylo_main_bound}
For every positive $\varepsilon\leq 1$ and a positive integer $r$,
there exist $\Delta_0$ and a constant $C=C(\varepsilon,r)$ such that:
$$\chi'_{a,r}(G)\leq\Delta(G)+C$$
for every graph $G$ without an isolated edge and with $\delta(G)\geq \varepsilon\Delta(G)$, $\Delta(G)\geq \Delta_0$.
In particular, $C\leq \varepsilon^{-2}(7r+200)+r+6$.
\end{theorem}

The general idea of proof of Theorem~\ref{przybylo_main_bound} was inspired by~\cite{Hatami}.
Its initial part is also similar to the second part of the proof of Theorem~\ref{first_attempt_theorem} above.

\section{Proof of Theorem~\ref{przybylo_main_bound}}

Let us fix $\varepsilon$ and $r$ as assumed, and let $G=(V,E)$ be a graph without isolated edges
and with maximum and minimum degrees $\Delta,\delta$, resp., such that $\delta\geq \varepsilon\Delta$.
Whenever needed we shall assume that $\Delta$ is sufficiently large, i.e., we explicitly do not
specify $\Delta_0$.
We shall prove that
$$\chi'_{a,r}(G)\leq\Delta+\varepsilon^{-2}(7r+200)+r+6.$$

Our randomized construction consists of two stages.\\
\\
\underline{\textbf{Stage One}:}\\
We first arbitrarily colour properly the edges of $G$ with
colours $1,2,\ldots,\Delta+1$. Denote this colouring by $c_0$. Then:
\begin{itemize}
\item we uncolour each edge $e\in E$ independently with probability $\frac{5r+100}{\varepsilon^2\Delta}$;
denote the set of uncoloured edges by $U$;
\item and finally we \emph{recover} every vertex with more than $\varepsilon^{-2}(7r+200)$ incident edges
uncoloured in the step above, i.e., we recover the removed colours of all the edges incident with such vertices.
\end{itemize}
Denote the partial colouring obtained by $c$, and let $U_c(v)$ denote the set of edges incident with $v\in V$ which are not coloured under $c$.

Note that $|U_c(v)|\leq \varepsilon^{-2}(7r+200)$ for every $v\in V$.
On the other hand, we shall prove that with positive probability every vertex $v$ with $|U_c(v)|\geq 3r+15$
is distinguished from all its $r$-neighbours (of the same degree) under $c$, while the vertices with $|U_c(v)|< 3r+15$ are rare and well distributed.
These shall be taken care of in stage two of the construction. Since we shall need some small additional alterations of $c$ within that stage,
we shall in fact prove below not only that (with positive probability) 
typically
the $r$-neighbours are distinguished after stage one, but also that this cannot change subject to later limited alterations of $c$, see
Claim~\ref{FinalClaimSt1} below.

Let $L$ be the set of all vertices $v\in V$ with $|U_c(v)|< 3r+15$.
Note that
\begin{equation}\label{Lincludedinsum}
L\subseteq R\cup LU\cup AR,
\end{equation}
where:

\begin{itemize}
\item $R$ is the set of all recovered vertices;
\item $LU$ is the set of vertices $v$ with $|U(v)|<3r+15$, where $U(v)$ is the set of edges incident with $v$ which belong to $U$;
\item $AR$ is the set of all vertices $v$ adjacent with some vertex $u\in R$ such that $uv\in U$.
\end{itemize}

The following technical observation shall be used several times in the further part of the argument.
\begin{observation}\label{two_technical_inequalities}
For every positive $\varepsilon \leq 1$ and a positive integer $r$,
\begin{eqnarray}
e^{10(1-\varepsilon^{-1})} &\leq& e^{10(1-\varepsilon^{-1})}\varepsilon^{-2} \leq \varepsilon^2,\label{EpsilonIneq}\\
2e^{-\frac{4r+170}{15}} &\leq& 2e^{-\frac{4r+170}{15}}(5r+100) \leq \frac{1}{100r}.\label{rIneq}
\end{eqnarray}
\end{observation}
\begin{pf}
Only the second inequality of (\ref{EpsilonIneq}) requires justification.
Note that it is equivalent to the fact that
$$f(\varepsilon)= 10(1-\varepsilon^{-1})-4\ln\varepsilon\leq 0,$$
which holds as
$$f'(\varepsilon)=\frac{4(2.5-\varepsilon)}{\varepsilon^2}$$
(hence $f$ is increasing for $\varepsilon\in(0,1]$) and $f(1)=0$.

Analogously, only the second inequality of (\ref{rIneq}) needs to be clarified.
Note that it is equivalent to the inequality:
\begin{equation}\label{EquivGIneq}
\frac{4r+170}{15} - \ln[(200r)(5r+100)]\geq 0,
\end{equation}
but since $(200r)(5r+100) = 1000r^2+20000r = (40r+160)^2 - 600(r-6)^2-4000 \leq (40r+160)^2$ (for $r\geq 1$),
to prove (\ref{EquivGIneq}) it is then sufficient to observe that
$$g(r)=\frac{4r+170}{15} - \ln(40r+160)^2\geq 0$$
for $r\geq 1$, as
$$g'(r)=\frac{4}{15}-\frac{2}{r+4}
=\frac{4(r-3.5)}{15(r+4)}$$
(i.e., $g$ is decreasing for $r\in[1,3.5]$ and
increasing for $r\geq 3.5$) and
$g(3.5)\approx 0.86 > 0$.
$\blacksquare$
\end{pf}

\begin{claim}\label{Claim_R}
For every vertex $v\in V$, ${\rm\emph{\textbf{Pr}}}\left(|N(v)\cap R| > \frac{\varepsilon^2\Delta}{30r} \right) <\frac{1}{30\Delta^{r+4}}$.
\end{claim}
\begin{pf}
Note first that by the Chernoff Bound, for any vertex $u\in V$,
\begin{eqnarray}
{\rm\textbf{Pr}}(u\in R)&=& {\rm\textbf{Pr}}(|U(u)|>\varepsilon^{-2}(7r+200))\nonumber\\
&\leq& {\rm\textbf{Pr}}\left({\rm BIN}\left(d(u),\frac{5r+100}{\varepsilon^2\Delta}\right)>\varepsilon^{-2}(7r+200)\right)\nonumber\\
&\leq& {\rm\textbf{Pr}}\left({\rm BIN}\left(\Delta,\frac{5r+100}{\varepsilon^2\Delta}\right)>\varepsilon^{-2}(7r+200)\right)\nonumber\\
&\leq& {\rm\textbf{Pr}}\left(\left|{\rm BIN}\left(\Delta,\frac{5r+100}{\varepsilon^2\Delta}\right)- \varepsilon^{-2}(5r+100)\right| > \varepsilon^{-2}(2r+100)\right)\nonumber\\
&\leq&2 e^{-\frac{\varepsilon^{-2}(2r+100)^2}{3(5r+100)}}\nonumber\\
&\leq& 2e^{-\frac{2\varepsilon^{-2}(2r+100)}{15}}\nonumber\\
&=& e^{-\frac{4r+200}{15}(\varepsilon^{-2}-1)}\cdot 2e^{-\frac{4r+200}{15}}\nonumber\\
&\leq& e^{10(1-\varepsilon^{-1})}\cdot 2e^{-\frac{4r+170}{15}}\nonumber\\
&\leq& \varepsilon^2\cdot\frac{1}{100r},\label{vinRineq}
\end{eqnarray}
where the last inequality follows by Observation~\ref{two_technical_inequalities}.

By (\ref{vinRineq}), for any vertex $v\in V$, ${\rm\textbf{E}}(|N(v)\cap R|)\leq \frac{\varepsilon^2\Delta}{100r}$.
We shall apply Talagrand's Inequality to the random variable $X=|N(v)\cap R|+\frac{\varepsilon^2\Delta}{100r}-{\rm\textbf{E}}(|N(v)\cap R|)$,
where ${\rm\textbf{E}}(X)=\frac{\varepsilon^2\Delta}{100r}$, to obtain the thesis.
For this aim, notice that $|N(v)\cap R|$ (and thus also $X$) is determined by the outcomes
of single trials associated with all edges incident with neighbours of $v$,
each of which sets down whether a colour is removed from
a given edge in the very first step of the construction or not.
Moreover, changing the outcome of each such trial may affect $|N(v)\cap R|$ (thus also $X$) by at most $2$,
and the fact that $|N(v)\cap R|\geq s$ (hence also that $X\geq s$) can be certified by the outcomes of at most
$(\varepsilon^{-2}(7r+200)+1)s$ trials. Therefore:
\begin{eqnarray*}
{\rm\textbf{Pr}}\left(|N(v)\cap R| > \frac{\varepsilon^2\Delta}{30r} \right) &\leq& {\rm\textbf{Pr}}\left(X > \frac{\varepsilon^2\Delta}{30r} \right)\\
&\leq& {\rm\textbf{Pr}}\left(X > \frac{\varepsilon^2\Delta}{100r}+\frac{\varepsilon^2\Delta}{100r}+120\sqrt{(\varepsilon^{-2}(7r+200)+1)\frac{\varepsilon^2\Delta}{100r}}\right)\\
&\leq& 4e^{-\frac{(\frac{\varepsilon^2\Delta}{100r})^2}{8\cdot 4\cdot (\varepsilon^{-2}(7r+200)+1)\frac{\varepsilon^2\Delta}{100r}}}\\
&<& \frac{1}{30\Delta^{r+4}}
\end{eqnarray*}
(for $\Delta$ sufficiently large).
$\blacksquare$
\end{pf}

\begin{claim}\label{Claim_LU}
For every vertex $v\in V$, ${\rm\emph{\textbf{Pr}}}\left(|N(v)\cap LU| > \frac{\varepsilon^2\Delta}{30r} \right) <\frac{1}{30\Delta^{r+4}}$.
\end{claim}
\begin{pf}
Note first that by the Chernoff Bound, for any vertex $u\in V$ ($d(u)\geq \varepsilon\Delta$),
for $\Delta$ sufficiently large,
\begin{eqnarray}
{\rm\textbf{Pr}}(u\in LU)&=& {\rm\textbf{Pr}}(|U(u)|<3r+15)\nonumber\\
&\leq& {\rm\textbf{Pr}}\left({\rm BIN}\left(\lceil\varepsilon\Delta\rceil,\frac{5r+100}{\varepsilon^2\Delta}\right)<3r+15\right)\nonumber\\
&\leq& {\rm\textbf{Pr}}\left(\left|{\rm BIN}\left(\lceil\varepsilon\Delta\rceil,\frac{5r+100}{\varepsilon^2\Delta}\right) - \lceil\varepsilon\Delta\rceil\frac{5r+100}{\varepsilon^2\Delta}\right| >  \varepsilon\Delta\frac{5r+100}{\varepsilon^2\Delta} - (3r+15)\right)\nonumber\\
&\leq& 2 e^{-\frac{[\varepsilon^{-1}(5r+100)-(3r+15)]^2}{3\lceil\varepsilon\Delta\rceil\cdot
\frac{5r+100}{\varepsilon^2\Delta}}}\nonumber\\
&\leq&2 e^{-\frac{[\varepsilon^{-1}(2r+85)]^2}{3\varepsilon^{-1}(5r+2.5\cdot 85)}}\nonumber\\
&=& e^{-\frac{4r+170}{15}(\varepsilon^{-1}-1)}\cdot 2e^{-\frac{4r+170}{15}}\nonumber\\
&\leq& e^{10(1-\varepsilon^{-1})}\cdot 2e^{-\frac{4r+170}{15}}\nonumber\\
&\leq& \varepsilon^2\cdot\frac{1}{100r},\label{vinLUineq}
\end{eqnarray}
where the last inequality follows by Observation~\ref{two_technical_inequalities}.

Now, instead of estimating the probability of $|N(v)\cap LU|$ being large,
we shall (equivalently) bound the probability that the random variable $X=d(v)-|N(v)\cap LU|$ is relatively small.
Note that by~(\ref{vinLUineq}),
$$\Delta \geq {\rm\textbf{E}}(X) = d(v)-{\rm\textbf{E}}(|N(v)\cap LU|) \geq d(v)-\frac{\varepsilon^2\Delta}{100r}.$$
Again, $X$
is determined by the outcomes
of single trials associated with all edges incident with neighbours of $v$,
each of which sets down whether a colour is removed from a given edge in the first step or not.
Moreover, changing the outcome of each such trial may affect $X$ by at most $2$,
and the fact that $X\geq s$ can be certified by the outcomes of at most
$(3r+15)s$ trials. Therefore, by Talagrand's Inequality:

\begin{eqnarray*}
&&{\rm\textbf{Pr}}\left(d(v)-|N(v)\cap LU|<d(v) - \frac{\varepsilon^2\Delta}{30r}\right)\\
&\leq&{\rm\textbf{Pr}}\left(X<\left(d(v)-\frac{\varepsilon^2\Delta}{100r}\right)-\frac{\varepsilon^2\Delta}{100r} - 120\sqrt{(3r+15)\Delta}\right)\\
&\leq&{\rm\textbf{Pr}}\left(|X-{\rm\textbf{E}}(X)|>\frac{\varepsilon^2\Delta}{100r} + 120\sqrt{(3r+15){\rm\textbf{E}}(X)}\right)\\
&\leq& 4e^{-\frac{(\frac{\varepsilon^2\Delta}{100r})^2}{8\cdot 2^2(3r+15)\mathbf{E}(X)}}\\
&=& e^{-\Omega(\Delta)}\leq \frac{1}{30\Delta^{r+4}}
\end{eqnarray*}
(for $\Delta$ sufficiently large).
The thesis follows.
$\blacksquare$
\end{pf}

\begin{claim}\label{Claim_AR}
For every vertex $v\in V$, ${\rm\emph{\textbf{Pr}}}\left(|N(v)\cap AR| > \frac{\varepsilon^2\Delta}{30r} \right) <\frac{1}{30\Delta^{r+4}}$.
\end{claim}
\begin{pf}
Note first that using the Chernoff Bound,
we obtain the following for any vertex $u\in V$ and its neighbour $w\in V$:
\begin{eqnarray}
{\rm\textbf{Pr}}(w\in R\wedge uw\in U)&=& {\rm\textbf{Pr}}(w\in R|uw\in U)\cdot {\rm\textbf{Pr}}(uw\in U)\nonumber\\
&\leq& {\rm\textbf{Pr}}\left({\rm BIN}\left(\Delta-1,\frac{5r+100}{\varepsilon^2\Delta}\right)>\varepsilon^{-2}(7r+200)-1\right)
\cdot \frac{\varepsilon^{-2}(5r+100)}{\Delta}\nonumber\\
&\leq& {\rm\textbf{Pr}}\left(\left|{\rm BIN}\left(\Delta,\frac{5r+100}{\varepsilon^2\Delta}\right)- \varepsilon^{-2}(5r+100)\right| > \varepsilon^{-2}(2r+99)\right) \cdot \frac{\varepsilon^{-2}(5r+100)}{\Delta}\nonumber\\
&\leq&2 e^{-\frac{\varepsilon^{-2}(2r+99)^2}{3(5r+100)}} \cdot \frac{\varepsilon^{-2}(5r+100)}{\Delta}.\nonumber
\end{eqnarray}
Therefore,
\begin{eqnarray}
{\rm\textbf{Pr}}(u\in AR)&\leq& \Delta\cdot 2 e^{-\frac{\varepsilon^{-2}(2r+99)^2}{3(5r+100)}} \cdot \frac{\varepsilon^{-2}(5r+100)}{\Delta}\nonumber\\
&\leq& 2e^{-\frac{2\varepsilon^{-2}(2r+99)}{15}} \cdot \varepsilon^{-2}(5r+100)\nonumber\\
&=& e^{-\frac{4r+198}{15}(\varepsilon^{-2}-1)} \varepsilon^{-2} \cdot 2e^{-\frac{4r+198}{15}} (5r+100)\nonumber\\
&\leq& e^{10(1-\varepsilon^{-1})}\varepsilon^{-2}\cdot 2e^{-\frac{4r+170}{15}}(5r+100)\nonumber\\
&\leq& \varepsilon^2\cdot\frac{1}{100r},
\label{vinARineq}
\end{eqnarray}
where the last inequality follows by Observation~\ref{two_technical_inequalities}.

By (\ref{vinARineq}), for any vertex $v\in V$, ${\rm\textbf{E}}(|N(v)\cap AR|)\leq \frac{\varepsilon^2\Delta}{100r}$.
We shall apply Talagrand's Inequality to the random variable $X=|N(v)\cap AR|+\frac{\varepsilon^2\Delta}{100r}-{\rm\textbf{E}}(|N(v)\cap AR|)$,
where ${\rm\textbf{E}}(X)=\frac{\varepsilon^2\Delta}{100r}$, to obtain the thesis.
For this aim, notice that $|N(v)\cap AR|$ (and thus also $X$) is determined by the outcomes
of single trials associated with all edges incident with neighbours of $v$ or their neighbours,
each of which
sets down whether a colour is removed from a given edge in the first step or not.
Moreover, changing the outcome of each such trial may affect $|N(v)\cap AR|$ (thus also $X$) by at most $2(\varepsilon^{-2}(7r+200)+1)$,
and the fact that $|N(v)\cap AR|\geq s$ (hence also that $X\geq s$) can be certified by the outcomes of at most
$(\varepsilon^{-2}(7r+200)+1)s$ trials. Therefore:
\begin{eqnarray*}
&&{\rm\textbf{Pr}}\left(|N(v)\cap AR| > \frac{\varepsilon^2\Delta}{30r} \right)\\
&\leq& {\rm\textbf{Pr}}\left(X > \frac{\varepsilon^2\Delta}{30r} \right)\\
&\leq& {\rm\textbf{Pr}}\left(X > \frac{\varepsilon^2\Delta}{100r}+\frac{\varepsilon^2\Delta}{100r} + 120(\varepsilon^{-2}(7r+200)+1) \sqrt{(\varepsilon^{-2}(7r+200)+1)\frac{\varepsilon^2\Delta}{100r}}\right)\\
&\leq& 4e^{-\frac{(\frac{\varepsilon^2\Delta}{100r})^2}{8\cdot 4\cdot (\varepsilon^{-2}(7r+200)+1)^3 \frac{\varepsilon^2\Delta}{100r}}}\\
&<& \frac{1}{30\Delta^{r+4}}
\end{eqnarray*}
(for $\Delta$ sufficiently large).
$\blacksquare$
\end{pf}

By (\ref{Lincludedinsum}) and Claims~\ref{Claim_R}, \ref{Claim_LU} and \ref{Claim_AR},
\begin{eqnarray}
{\rm\textbf{Pr}}\left(|N(v)\cap L| > \frac{\varepsilon^2\Delta}{10r}\right) &\leq&  {\rm\textbf{Pr}}\left(|N(v)\cap R| > \frac{\varepsilon^2\Delta}{30r} \vee |N(v)\cap LU| > \frac{\varepsilon^2\Delta}{30r} \vee |N(v)\cap AR| > \frac{\varepsilon^2\Delta}{30r}\right) \nonumber\\
&<& \frac{1}{10\Delta^{r+4}}\label{NvL_bound}
\end{eqnarray}
for every $v\in V$.

Let $A \triangle B$ denote the symmetric difference of any two sets $A$ and $B$, i.e., $A\triangle B=(A\smallsetminus B)\cup(B\smallsetminus A)$.

\begin{claim}\label{Claim_symm_diff}
For every $u,v\in V$ with $d(u)=d(v)$ and $1\leq d(u,v)\leq r$,
$${\rm\emph{\textbf{Pr}}}(u\notin L \wedge |S_c(u)\triangle S_c(v)|<2r+10)<\frac{1}{10\Delta^{r+4}}.$$
\end{claim}

\begin{pf}
Consider any vertices $u,v\in V$ with $d(u)=d(v)$ and $1\leq d(u,v)\leq r$.
We wish to upper-bound the probability of the event: $u\notin L \wedge |S_c(u)\triangle S_c(v)|<2r+10$. 
Note then that $u\notin L$ in particular implies that $3r+15\leq |U(u)|\leq \varepsilon^{-2}(7r+200)$ and $|U_c(v)|\geq 3r+15$.
Since our random uncolourings of the edges are independent, we may consider these being done in any order. 
Suppose we first perform the corresponding experiments (determining whether a given edge is uncoloured or not)
for the edges incident with $u$. 
%
%
Afterwards, there are at least $|U(u)|$ (or $|U(u)|-1$ if $uv\in E$ and $uv$ was uncoloured) edges incident with $v$ 
whose colours do not belong to the pallet of $u$. Among these choose $|U(u)|$ ($|U(u)|-1$, resp.) with the least colours 
(recall that we use integer colours in our construction) and denote them by $E'_v$.
Since we must have that $|U_c(v)|\geq 3r+15$ within the investigated event, at most $|U(u)|-3r-15$ edges in $U(u)$ might have their 
colours recovered 
eventually (some or all of which might have been assigned to the edges in $E'_v$). 
In order to have $|S_c(u)\triangle S_c(v)|<2r+10$, i.e., $|S_c(u)\triangle S_c(v)|\leq 2r+9$ at the end,
still at least $r+5$ edges in $E'_v$ must be uncoloured in our random process. 
Since we also must have $|U(u)| 
\leq \varepsilon^{-2}(7r+200)$, hence $|E'_v| \leq \varepsilon^{-2}(7r+200)$, we obtain that:
\begin{eqnarray}
{\rm\textbf{Pr}}(u\notin L \wedge |S_c(u)\triangle S_c(v)|<2r+10) 
&\leq& {\lfloor\varepsilon^{-2}(7r+200)\rfloor \choose r+5} \left(\frac{\varepsilon^{-2}(5r+100)}{\Delta}\right)^{r+5} \nonumber\\
&<& \frac{1}{10\Delta^{r+4}} \nonumber
\end{eqnarray}
(for $\Delta$ sufficiently large).
$\blacksquare$
\end{pf}

\begin{claim}\label{FinalClaimSt1}
We can choose the partial colouring $c$ so that:
\begin{itemize}
\item[(a)] $|N(v)\cap L| \leq \frac{\varepsilon^2\Delta}{10r}$
for every vertex $v\in V$, and
\item[(b)] $|S_c(u)\triangle S_c(v)|\geq 2r+10$ for every pair of $r$-neighbours $u,v$
with $d(u)=d(v)$ and $u\notin L$ (or $v\notin L$).
\end{itemize}
\end{claim}

\begin{pf}
For every $u,v\in V$, let $A_v$ denote the event that $|N(v)\cap L| > \frac{\varepsilon^2\Delta}{10r}$,
and if $d(u)=d(v)$ and $1\leq d(u,v)\leq r$, let $A_{u,v}$ denote the event
that $u\notin L$ and $|S_c(u)\triangle S_c(v)|<2r+10$
(note that $A_{u,v}$ differs from $A_{v,u}$ within this convention).
Observe that every event $A_v$ is mutually independent of all events $A_{v'}$ and $A_{u,w}$
with $d(v,v')>5$, $d(v,u)>4$ and $d(v,w)>4$, i.e., of all other events of these forms
but at most $\Delta^5+\Delta^4\cdot\Delta^r\cdot 2 \leq 3\Delta^{r+4}$.
Analogously, every event $A_{u,v}$ is mutually independent of all events $A_w$ and $A_{u',v'}$
with $d(u,w)>4$, $d(v,w)>4$ and $d(u,u')>3$, $d(u,v')>3$, $d(v,u')>3$, $d(v,v')>3$,
i.e., of all other events of these forms
but at most $2\Delta^4+2\cdot\Delta^3\cdot\Delta^r\cdot 2 \leq 3\Delta^{r+4}$.
Moreover, by~(\ref{NvL_bound}) and Claim~\ref{Claim_symm_diff},
each of these events occurs with probability at most $\frac{1}{10\Delta^{r+4}}$.
By the Lov\'asz Local Lemma,
we may thus perform the uncolourings so that none of these events holds for the obtained partial colouring $c$.
$\blacksquare$
\end{pf}
\underline{\textbf{Stage Two}:}\\
Note that after stage one, all vertices not in $L$ are distinguished from their $r$-neighbours (of the same degrees).
Now we shall slightly modify our colouring to
ensure the same for the vertices in $L$.
Thus for every $v\in L$ we randomly uncolour its $r+5$
incident edges which were coloured under $c$ (obtained after stage one) joining $v$ with vertices outside $L$.
The choices are independent for all vertices in $L$ (and feasible, as
$d(v)-|U_c(v)|-|N(v)\cap L| > \varepsilon\Delta - (3r+15) - \frac{\varepsilon^2\Delta}{10r}$
for every $v\in L$).
We denote the partial colouring obtained by $c'$, the set of uncoloured within this stage edges by $U'$,
and the subset of these incident with any given $v\in V$ by $U'(v)$.
Obviously such uncolourings
influence also the colour pallets of vertices outside $L$,
but we shall show that these
changes may be minor (with positive probability),
and the condition
\emph{(b)} of Claim~\ref{FinalClaimSt1} above suffices to keep all these vertices distinguished from their $r$-neighbours (with the same degrees).
\begin{claim}\label{mainSecondStageClaim}
We can perform the changes in the colouring $c$ described above so that for the obtained partial colouring $c'$:
\begin{itemize}
\item[(a)] for every vertex $u\in V\smallsetminus L$,  $|U'(u)|\leq r+4$;
\item[(b)] $S_{c'}(u)\neq S_{c'}(v)$ for every $u,v\in L$ with $d(u)=d(v)$ and $1\leq d(u,v)\leq r$.
\end{itemize}
\end{claim}
\begin{pf}
We define two kinds of bad events:
\begin{itemize}
\item for every vertex $u\in V\smallsetminus L$
  and any its $r+5$ neighbours $v_1,\ldots,v_{r+5}\in L$ such that
  $uv_i\notin U_c(u)$ for $i=1,\ldots,r+5$,
let $A_{u,\{v_1,\ldots,v_{r+5}\}}$ denote the event that $uv_i\in U'$ for $i=1,\ldots,r+5$;
\item for every $u,v\in L$ with $d(u)=d(v)$ and $1\leq d(u,v)\leq r$, let $B_{\{u,v\}}$ denote the event that $S_{c'}(u)=S_{c'}(v)$.
\end{itemize}

Since by stage one, every vertex $v\in L$ has at least
$\varepsilon\Delta-\frac{\varepsilon^{2}\Delta}{10r}-(3r+15)$ neighbours $w\notin L$
such that
$vw\notin U_c(v)$, then in all cases:
\begin{equation}\label{Auvi_ineq}
{\rm\textbf{Pr}}(A_{u,\{v_1,\ldots,v_{r+5}\}})\leq  \left(\frac{r+5}{\varepsilon\Delta-\frac{\varepsilon^{2}\Delta}{10r}-(3r+15)}\right)^{r+5}\leq
\left(\frac{r+5}{\varepsilon\Delta-\frac{2\varepsilon\Delta}{10r}}\right)^{r+5}
\leq \left(\frac{r+5}{\frac{4}{5}\varepsilon\Delta}\right)^{r+5},
\end{equation}
\begin{equation}\label{Buv_ineq}
{\rm\textbf{Pr}}(B_{\{u,v\}}) \leq \frac{1}{{\lceil\varepsilon\Delta-\frac{\varepsilon^{2}\Delta}{10r}-(3r+15)\rceil \choose r+5}}
\leq \frac{1}{\left(\frac{\frac{4}{5}\varepsilon\Delta}{r+5}\right)^{r+5}}
= \left(\frac{r+5}{\frac{4}{5}\varepsilon\Delta}\right)^{r+5}.
\end{equation}

Let us define a dependency graph $G'$ with vertex set consisting of all the events of both types above as follows.
Every event $A_{u,\{v_1,\ldots,v_{r+5}\}}$ is adjacent with $A_{u',\{v'_1,\ldots,v'_{r+5}\}}$ if $\{v_1,\ldots,v_{r+5}\}\cap\{v'_1,\ldots,v'_{r+5}\}\neq \emptyset$, while it is adjacent with $B_{\{u',v'\}}$ if $\{v_1,\ldots,v_{r+5}\}\cap\{u',v'\}\neq \emptyset$, and analogously, $B_{\{u,v\}}$ and $B_{\{u',v'\}}$ form an edge in $G'$ if and only if $\{u,v\}\cap\{u',v'\}\neq\emptyset$.
Then every $A_{u,\{v_1,\ldots,v_{r+5}\}}$, and analogously each $B_{\{u,v\}}$ is mutually independent
of all other events which are not adjacent with it in $G'$.
In order to upper-bound the degrees of these events in $G'$, recall that by Claim~\ref{FinalClaimSt1}\emph{(a)},
each vertex $v$ of $G$ has at most $\frac{\varepsilon^2\Delta}{10r}$ neighbours in $L$.
Obviously, it can also have at most $\Delta^r$ $r$-neighbours in $G$ (which might belong to $L$).
Therefore, for all the investigated events, we have
$d_{G'}(A_{u,\{v_1,\ldots,v_{r+5}\}}) \leq (r+5)\Delta {\lfloor\frac{\varepsilon^2\Delta}{10r}\rfloor \choose r+4} + (r+5)\Delta^{r}$ and
$d_{G'}(B_{\{u,v\}}) \leq 2\Delta {\lfloor\frac{\varepsilon^2\Delta}{10r}\rfloor \choose r+4} + 2\Delta^{r}$,
hence
\begin{equation}\label{DependencyDeltaBound}
\Delta(G') \leq (r+5)\left(\Delta {\lfloor\frac{\varepsilon^2\Delta}{10r}\rfloor \choose r+4} + \Delta^{r}\right)
\leq (r+6) \Delta {\lfloor\frac{\varepsilon^2\Delta}{10r}\rfloor \choose r+4}
\leq \frac{r+6}{(r+4)!}\Delta\left(\frac{\varepsilon^2\Delta}{10r}\right)^{r+4}.
\end{equation}

As
$$e \left(\frac{r+5}{\frac{4}{5}\varepsilon\Delta}\right)^{r+5} (\Delta(G')+1)
< 4\left(\frac{r+5}{\frac{4}{5}\varepsilon\Delta}\right)^{r+5}
\frac{r+6}{(r+4)!}\Delta\left(\frac{\varepsilon^2\Delta}{10r}\right)^{r+4}
\leq \frac{5(r+6)(r+5)^{r+5}}{(r+4)!(8r)^{r+4}}
< 1$$
(where the last inequality can be checked directly for $r=1$, while for $r\geq 2$ we obviously have:
$(r+5)^{r+4}<(8r)^{r+4}$ and $5(r+6)(r+5)<(r+4)!$),
by (\ref{Auvi_ineq}), (\ref{Buv_ineq}), (\ref{DependencyDeltaBound}) and the Lov\'asz Local Lemma,
we may choose $c'$ so that none of the events of the forms
$A_{u,\{v_1,\ldots,v_{r+5}\}}$ and  $B_{\{u,v\}}$
holds.
$\blacksquare$
\end{pf}

By our construction, all $r$-neighbours of the same degrees are distinguished under $c'$, cf. Claims~\ref{FinalClaimSt1} and~\ref{mainSecondStageClaim}.
Additionally, for every vertex $v\notin L$, $d(v)-d_{c'}(v) \leq \varepsilon^{-2}(7r+200) + (r+4)$,
while if $v\in L$, then $d(v)-d_{c'}(v) < (3r+15) + (r+5)\leq \varepsilon^{-2}(7r+200) + (r+4)$.
By Vizing's Theorem, we may use at most $\varepsilon^{-2}(7r+200) + (r+4)+1$ new colours
for the uncoloured edges
in order to extend the colouring $c'$ to a proper colouring of the whole graph $G$. Counting in the initial (at most) $\Delta+1$ colours,
we have thus used in total no more than $\Delta+\varepsilon^{-2}(7r+200) + (r+4)+2$ colours
to construct an $r$-distant set distinguishing colouring of $G$.
The proof of Theorem~\ref{przybylo_main_bound} is thus completed.
$\blacksquare$

\section{Comments}

Note that
substituting $\varepsilon = 1$ in Theorem~\ref{przybylo_main_bound} we obtain
the following.
\begin{corollary}\label{regular_graphs_corollary}
For every positive integer $r$, there exists $d_0$ such that:
$$\chi'_{a,r}(G)\leq \Delta(G)+8r+206$$
for every $d$-regular graph $G$ with $d\geq d_0$.
\end{corollary}
This in particular implies that Conjecture~\ref{Przybylo_main_conjecture} holds for regular graphs.
This also proves that for each fixed $r$, $\chi'_{a,r}(G)\leq \Delta(G)+C$
for every regular graph without isolated edges, where $C$ is some constant dependent on $r$.
For large degrees it follows by Corollary~\ref{regular_graphs_corollary} above.
In the remaining cases, i.e., for any graph $G$ with $\Delta(G) < d_0$
one can easily prove that $\chi'_{a,r}(G)\leq C_0$, where $C_0=C_0(d_0,r)$ is some (large enough) constant.
It is sufficient to use a greedy approach exploiting
induction based on removal of two incident
edges from a given graph.
Analogously, by Theorem~\ref{przybylo_main_bound}, Conjecture~\ref{Przybylo_main_conjecture} holds
for isolated edge free graphs $G$ with $\delta(G)\geq \varepsilon\Delta(G)$,
where $\varepsilon\leq 1$ is any fixed positive constant.
By a similar greedy argument as above (and  Theorem~\ref{przybylo_main_bound}),
$\chi'_{a,r}(G)\leq \Delta(G)+C$ for all such graphs,
where $C$ is some constant dependent on any fixed $r$ and $\varepsilon$.

It would be interesting to develop such greedy approach
towards designing a general upper bound on $\chi'_{a,r}$ independent of the minimum degree,
improving the known bounds from~\cite{Tian}.
Also the quest for
an upper bound of the form $\chi'_{a,r}(G) \leq \Delta(G)+const.$ (for each fixed $r$)
for all graphs with minimum degree larger than a constant dependent only on $r$
remains an open problem in general.

At the end it is also worth mentioning that similar results can be achieved in the case of total colourings, see~\cite{Przybylo_distant_set_total}.


\end{document}